\documentclass[11pt,twoside]{article}
\usepackage{amsmath, amsthm, amsfonts, amssymb}
\usepackage[all]{xy}
\def\email#1{\vspace{1cm}\noindent E-mail address: {\sf #1} \bigskip}

\theoremstyle{plain}
\newtheorem{thm}{Theorem}
\newtheorem{proposition}[thm]{Proposition}
\newtheorem{corollary}[thm]{Corollary}
\newtheorem{lemma}[thm]{Lemma}
\theoremstyle{remark}
\newtheorem{remark}[thm]{Remark}

\theoremstyle{definition}
\newtheorem{definition}[thm]{Definition}
\newtheorem{example}[thm]{Example}
\newtheorem{exap}[thm]{Example-Application}

\newtheorem{defirema}[thm]{Definition and Remark}

\def\nerve{\mathcal N}

\def\set#1{\{#1\}}
\def\u{\mathcal U}
\def\w{\mathcal V}
\def\yr{(Y,R)}
\def\zr{(Z,\tilde R)}
\def\k{K_*}
\def\ky{K_Y}
\def\ly{L_Y}
\def\kz{K_Z}
\def\lz{L_Z}

\title{The geometry of relations}

\author{Elias Gabriel Minian}

\date{}

\begin{document}

\maketitle

{Departamento  de Matem\'atica.\\
 FCEyN, Universidad de Buenos Aires. \\ Buenos
Aires, Argentina}

\begin{abstract}
\noindent 
There is a canonical way to associate two simplicial complexes $K, L$ to any relation $R\subset X\times Y$. Moreover, the geometric
realizations of $K$ and $L$  are  homotopy
equivalent.
This was studied in the fifties by C.H. Dowker \cite{d}. In this article we
prove a Galois-type correspondence for relations $R\subset X\times Y$ when $X$ is fixed and  use these constructions 
to investigate finite posets (or equivalently, finite $T_0$-spaces) from a geometrical point of view. Given any poset $(X, \leq)$, 
we define the  simplicial complexes $K, L$ 
associated to the 
relation $\leq$. In 
many cases these polyhedra have the same homotopy type as the standard 
simplicial complex $C_X$ of nonempty finite chains in $X$. We give a complete characterization of the simplicial complexes 
that are the $K$ or $L$-complexes of
some finite poset and prove that $K$ and $L$ are geometrically equivalent 
to the smaller complexes
$K',L'$ induced by the relation $<$. More precisely, we  
prove that $K$ (resp. $L$) simplicially collapses to $K'$ (resp. $L'$). 
\end{abstract}



\bigskip

\bigskip

Let $R$ be a relation between two nonempty sets $X$ and $Y$, i.e. $R$ is a subset of the cartesian product $X\times Y$. We write $xRy$ if 
$(x,y)\in R$. There is a canonical way to associate to $R$ two simplicial complexes $K$ and $L$. The simplicial complex $K$ is defined as follows.
The $n$-simplices of $K$ are the finite subsets $\{x_0,\ldots,x_n\}$ of $X$ such that there exists some $y\in Y$ with $x_iRy$ for all $i=0,\ldots,n$.
Similary, the simplices of  $L$ are the finite subsets $s$ of $Y$ such that the elements of $s$ are related to a common element of $X$. In particular, 
the set of vertices of $K$ is a subset of $X$ and consists of the points of $X$ which are related to some element in $Y$. Analogously, the set
of  vertices of $L$ is a subset of $Y$. Sometimes it is useful to refer to $K$ and $L$ as the $K$-complex and the $L$-complex associated to the
relation $R$.

\bigskip

These constructions were introduced by C.H. Dowker in \cite{d} where he proved the following result.

\begin{thm}[C.H. Dowker]
Let $R\subset X\times Y$ be a relation and let $K$ and $L$ be the associated simplicial complexes. Then the polyhedra $|K|$ and $|L|$ are homotopy
equivalent.
\end{thm}

Here $|K|$ denotes the geometric realization of $K$.

\bigskip

Dowker applied this result to relate the \v{C}ech and Vietoris homology groups of a space. More precisely, let $X$ be a topological space and $\u$
a cover of $X$, i.e. $X=\bigcup_{U\in \u}U$. Consider the relation $xR U$ if $x\in U$. The $K$-complex is called in this case the Vietoris
complex of the covering, and denoted $V(\u)$, and the $L$-complex  is called the nerve of the covering, which is  denoted by $\nerve(\u)$. 
Since the geometric realizations
of $V(\u)$ and $\nerve(\u)$ are homotopy equivalent, it follows that the \v{C}ech and Vietoris homology coincide.

\bigskip

Another interesting applications are given in \cite{gaga}. In that paper these constructions are used 
 to study  the homotopy theory of \it groupoid atlases. \rm

\bigskip

The first goal of this article is to investigate the relationship between the different relations $R\subset X\times Y$ 
defined on a fixed nonempty set $X$. We are interested
 to study for example the relationship between different coverings of a given space $X$,  or the different  poset structures that can be defined 
on a given set. We will prove a Galois-type correspondence between the subcomplexes of the simplicial complex spanned by the set $X$ and the equivalence 
classes of relations $R\subset X\times Y$.

\medskip

In the second part of the paper we will investigate finite posets (or equivalently, finite $T_0$-spaces) from a geometrical viewpoint using various 
simplicial complexes associated to  the relations $\leq$ and $<$ of the poset.

\bigskip

We introduce first some notations. 

Let $X$ be a fixed nonempty set. We denote by $\yr$ a relation $R\subset X\times Y$. Let $\ky$ be the $K$-complex and $\ly$ the $L$-complex 
 associated
to $\yr$. We will work with \it covered \rm relations $\yr$, i.e. relations such that the proyection on $Y$, $p:R\to Y$ is onto or equivalently, 
for any $y\in Y$ there is an element $x\in X$ such that $xRy$.

\medskip
We denote by $\k$ the simplicial complex spanned by the set $X$, i.e. the simplices of $\k$ are all finite subsets of $X$. 
When $T$ is a subcomplex of some simplicial complex $K$, we will write $T\leq K$.
Note that
for any $\yr$, $\ky\leq \k$. Conversely, we have

\begin{proposition}\label{asignarelacion}
For any  $T\leq\k$, there exists a (covered) relation $\yr$ such that $\ky=T$.
\end{proposition}

\begin{proof}

Take $Y=S_T$ the set of simplices of $T$ and define $xRs$ if $x\in s$. It is easy to verify that $\ky=T$
\end{proof}

\begin{definition}
A morphism $f:\yr\to\zr$ is a set theoretic map $f:Y\to Z$ such that for every $x\in X$ and $y\in Y$
$$xRy\Rightarrow x\tilde R f(y)$$
\end{definition}

\begin{remark}\label{basico}
Note that a morphism $f:\yr\to\zr$ induces a well defined simplicial map $L_f:\ly\to\lz$. Moreover, if there exists a morphism
$f:\yr\to\zr$, then $\ky$ is a subcomplex of $\kz$.
\end{remark}

\begin{proposition}\label{contiguous}
Let $f,g:\yr\to\zr$ be morphisms. Then the induced simplicial maps $L_f,L_g:\ly\to \lz$ are contiguous. In particular, $|L_f|$ and 
$|L_g|$ are homotopic continuous maps.
\end{proposition}

\begin{proof}
It is sufficient to verify that for any simplex $s=\set{y_0,\ldots,y_n}$ of $\ly$, $L_f(s)\cup L_g(s)$ is a simplex in $\lz$. 

Since $s$ is a simplex, there exists $x\in X$ such that $xRy_j$ for every $j$ and therefore $x\tilde Rf(y_j)$ and
$x\tilde R g(y_j)$ for every $j$. The result follows immediately.
\end{proof}

\begin{definition}
Two relations $\yr$ and $\zr$ are called \it equivalent \rm if there are morphisms $f:\yr\to\zr$ and $h:\zr\to \yr$. We denote in this case
$\yr\simeq \zr$.
\end{definition}

\begin{remark}
Note that $\simeq$ is an equivalence relation in the class of relations defined on $X$. To prove that it is transitive, note that any composition of
morphisms is again a morphism.
\end{remark}

From remark \ref{basico} and proposition \ref{contiguous}, it follows straightforward the following

\begin{corollary}\label{equivalent}
If $\yr\simeq \zr$, then
\begin{enumerate}
\item $\ky=\kz$.
\item  $|\ly|$ and $|\lz|$ are homotopy equivalent.

\item Any morphism $g:\yr\to\zr$ induces a homotopy equivalence $|L_g|$.
\end{enumerate}
\end{corollary}

Note that, if the hypothesis $\yr\simeq \zr$ is not satisfied, a morphism $g:\yr\to\zr$ does not in general induce a homotopy equivalence. Consider for
example the constant morphism $\yr\to (\ast,\tilde R)$ where $\ast$ is the singleton and $x\tilde R \ast$ for all $x\in X$. 

\medskip

Note also that if $g:\yr\to\zr$ induces a homotopy equivalence, this does not imply that $\yr$ and $\zr$ are equivalent. Consider for instance a relation
$\yr$ such that $\ly$ is contractible. In that case, the map $\yr\to (\ast,\tilde R)$ induces a homotopy equivalence but in general $\yr$ is not
equivalent to $(\ast,\tilde R)$ unless there exists $y\in Y$ such that $xRy$ for all $x\in X$.

\begin{defirema}
Given a relation $\yr$ and an element $y\in Y$, let $S_y$ be the set of all elements of $X$ which are related to $y$. Since we work 
with covered relations,
these sets are nonempty. Moreover $S_y$ is a generalized simplex of $\ky$, i.e. all its finite subsets are simplices of $\ky$. If $S_y$ is finite, it is 
just a simplex in $\ky$.
\end{defirema}

We prove now the Galois-type correspondence for relations:

\begin{thm}
Let $X$ be a finite set and let $K_*$ be the simplicial complex spanned by $X$. There exists a one-to-one correspondence between subcomplexes of
$K_*$ and equivalence classes of covered relations $\yr$ on $X$. This correspondence assigns to each subcomplex $T$ of $K_*$ the class 
of the relation $(S_T,R)$ as
in proposition \ref{asignarelacion} and to each class $\yr$ the subcomplex $\ky\leq \k$. Moreover, 
$$\ky\leq \kz\Longleftrightarrow \exists\  f:\yr\to\zr$$
 \end{thm}

\begin{proof}
The first part of the Theorem follows from results \ref{asignarelacion}, \ref{basico} and \ref{equivalent} of above. 

\medskip

Suppose now that $\yr$ and $\zr$ are relations on $X$ with $\ky\leq\kz$. Let $y\in Y$ and let $S_y$ be its associated generalized simplex. Since 
$X$ is finite and $\ky\leq\kz$, then $S_y$ is actually a simplex of $\kz$ and therefore there exists some $z\in Z$ such that $x\tilde R z$ for all
$x\in S_y$. Define $f(y)=z$. The function $f:Y\to Z$ defined this way is a morphism of relations since
$$xRy\Longrightarrow x\in S_y \Longrightarrow x\tilde R f(y).$$
This completes the proof.
\end{proof}

Note that the finiteness hypothesis is needed to prove the implication 
$$\ky\leq\kz \Longrightarrow \exists\  f:\yr\to\zr.$$
If $X$ is infinite, take $Y=P(X)$ the set of all nonempy subsets of $X$ and let $Z=P_f(X)$ the set of all nonempty finite subsets of $X$ and define
$xRT$ (or $x\tilde RT$) if $x\in T$. Then $\ky=\kz$ but there is no morphism $f:(Y,R)\to\zr$.

\medskip

\bigskip

{\bf Posets, finite spaces and nerves.}

\bigskip

We exhibit now some examples and applications.

\begin{exap}\label{refines}
Take a topological space $X$ and two covers $\u$ and $\w$ and suppose that $\w$ refines $\u$ (i.e. for any $V\in\w$ there is some $U\in \u$ such that
$V\subseteq U$). If we consider $\u$ and $\w$ as relations on $X$ (as above), then there is a morphism of relations $f:\w\to\u$. In fact such a morphism
is precisely a \it refinement map \rm (or a \it canonical projection \rm  in the terminology of \cite{sp}). Therefore  $\w$ refines $\u$ if and only if there is a morphism
of relations $f:\w\to\u$ and two coverings refine each other if and only if they are equivalent (as relations on $X$). As an immediate corollary of 
proposition 
\ref{equivalent}, we verify the 
well-known fact that the nerves of  two covers which refine each other are homotopy equivalent.

\end{exap}

\bigskip

Consider now a finite poset (=partially ordered set) $(X,\leq)$ and the relation $R\subset X\times X$ given by $\leq$. As usual, denote by
  $K$ 
and $L$ the complexes associated to $\leq$. Note that the simplices of $K$ are the subsets $\set{x_0,\ldots,x_n}$ of $X$ such that there exists $y\in X$ with $x_i\leq y$ for all $i$. 
Similarly, the simplices of $L$ are the subsets with a common lower bound $z\in X$.

\bigskip

Sometimes it is very useful to regard finite posets  as finite $T_0$-spaces.
More precisely, given a finite poset $(X,\leq)$, we define for each $x\in X$ the 
set $$U_x=\{y\in X \ | \ y\le x\}.$$
It is not difficult to verify that these sets form a basis for a topology. This is the topology associated to $\le$.

Conversely, given  a  finite $T_0$ topological space $X$, we consider for each  $x \in X$ the \textit{minimal open set} $U_x$ which
 is defined as the intersection of all open sets containing $x$. The partial order associated to the topology on $X$ is given by 
the relation $x\le y$ if $x\in U_y$.

These applications define a one-to-one 
correspondence between $T_0$-topologies and partial orders on the finite set $X$. Therefore one can consider 
finite $T_0$-spaces as finite posets and viceversa.

For more details on finite spaces, some references are \cite{bm, bm2, M, M2, Mcc1, Sto}.

\bigskip

The poset $(X,\leq)$ can be identified with the poset $(\u,\subseteq)$ of minimal open sets of $X$ (ordered by inclusion), each element $x$ is identified
with its minimal open set $U_x$. Analogously, $X$ can be identified with the poset $(\w,\supseteq)$ of minimal closed sets of $X$. Explicitly, 
$$\w=\set{F_x,\ x\in X}$$
 where $F_x$ is the intersection of all closed subsets of $X$ containing $x$.

Under these identifications, the simplicial complexes $K$ and $L$ associated to the relation $\leq$ defined on $X\times X$, are exactly the nerves
$\nerve(\w)$ and $\nerve(\u)$ of the minimal closed and open sets respectively. By Dowker's Theorem, we obtain

\begin{corollary}
Let $X$ be a finite $T_0$-space and let $\nerve(\w)$ and $\nerve(\u)$ denote the nerves of the minimal closed and open sets of $X$. Then
 $$| \nerve(\w) | \simeq | \nerve(\u) |.$$
\end{corollary}

\medskip

The rest of the paper is devoted to investigate the information that one can obtain from the simplicial complexes $K$ and $L$ associated to 
a finite poset
$(X,\leq)$. 

\bigskip

There is a classical way to associate a simplicial complex $C_X$ to a (finite) poset $X$ (see for example \cite{Bjo0,Bjo, Bjo2, Qui}): The
simplices of $C_X$ are the nonempty finite chains in $X$. Moreover, if we regard $X$ as a finite space, there is a weak homotopy equivalence 
$|C_X|\to X$ (see \cite{bm, M2, Mcc1}). By definition it is clear that $C_X$ is a subcomplex of $K$ and $L$ but in general it does not have the
same homotopy type of $K$ (and $L$). Consider for example the finite poset $X$ with Hasse diagram

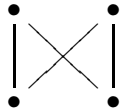
\begin{displaymath}
\xymatrix@C=10pt{ \bullet \ar@{-}[d] \ar@{-}[drr] & & \bullet \ar@{-}[lld]  \ar@{-}[d]  \\
		\bullet & & \bullet  } 
\end{displaymath}

Clearly $C_X$ is a one dimensional sphere $S^1$ but $K$ and $L$ are contractible.

\medskip

In many situations, $C_X$ is homotopy equivalent to $K$ and $L$, moreover the inclusions $C_X\rightarrowtail K$ and $C_X\rightarrowtail L$ are deformation retracts. This 
follows from the following  result of McCord \cite{Mcc2}:

\begin{lemma}[McCord]
Let $\u$ be a cover of a space with the property that the intersection of any finite collection of elements of $\u$ is either empty or a member of 
$\u$. Let $C(\u)$ be the subcomplex of the nerve $\nerve(\u)$ whose simplices are the nonempty chains of $\u$. Then the inclusion 
$C(\u)\rightarrowtail \nerve(\u)$ is a deformation retract.
\end{lemma}
 
If we take  $\u$ as the covering of minimal open sets of $X$ and use the identifications of above, we obtain

\begin{corollary}\label{rl}
Let $X$ be a finite $T_0$-space such that for any $x,y\in X$ the intersection $U_x\cap U_y$ is either empty or $U_z$ for some $z\in X$. Then the inclusion
$C_X\subset L$ is a deformation retract.
\end{corollary}

Therefore for posets $X=\mathcal{L}-\set{\hat 0, \hat 1}$, with $\mathcal{L}$ a finite lattice, 
 the simplicial complexes $K$ and $L$ have the same homotopy type as the standard simplicial complex $C_X$. Note that 
Corollary \ref{rl} follows
also from the Crosscut Theorem (see for example \cite{Bjo0}, Thm. 10.8) and Example-Application \ref{refines}.

\bigskip

Suppose now that $X$ is a finite set. We want to relate the constructions of above with the classification of poset 
structures that can be defined on the set $X$. 
More explicitly: Given a subcomplex $T$ of the simplicial complex $\k$ spanned by $X$, can we define a poset structure $\leq$ on $X$ such that $T$ is the
$K$-complex (resp. the $L$-complex) of this relation?

\medskip

In order to be the $K$-complex of a poset structure $\leq$ on $X$, $T$ should be, first of all, \it complete, \rm i.e. the set of vertices of $T$ should be
the whole set $X$ since the relation $\leq$ must be reflexive. The second condition that $T$ must satisfy is deduced from the following lemma.

\begin{lemma}\label{maximal}
Let $X$ be a finite poset and let $T$ be the associated
$K$-complex. Let $s$ be a  maximal simplex of $T$. Then $s=\set{y}$ or $s=\set{y,x_0,\ldots,x_n}$ where $y$ is a maximal element in $X$ and 
the (possibly empty)
set $\set{x_0,\ldots,x_n}$ consists  of all elements of $X$ such that $x_i < y$ (equivalently, the set of all elements of $X$ which
are comparable with $y$). In particular, any maximal simplex $s$ of $T$ contains exactly one maximal element of $X$ and any maximal element of $X$ is
in only one maximal simplex of $T$.
\end{lemma}

\begin{proof}
Let $s=\set{z_0,\ldots,z_m}$ be a maximal simplex of $T$. Since $s\in T$, then there exists some $y\in X$ such that $z_i\leq y$ for all $i$. Since $s$ 
is maximal, $y=z_j$ for some $j$ and $y$ must be also a maximal element of $X$, for if $y< w$ for some $w$, then we could add $w$ to $s$ and this 
contradicts
the maximality of $s$. Therefore one (and only one) of the $z_i's$ is a maximal element $y$ of $X$ and the others are the elements of $X$ which are 
smaller than $y$ (since $s$ is maximal,  all elements smaller than $y$ must belong to $s$).
\end{proof}

Thus, in order to be the $K$-complex of a poset structure, in any maximal simplex $s$ of $T$ there must be an element $y$ which does not belong to
any other maximal simplex of $T$. For example, the boundary of the closed 2-simplex is a simplicial complex of 3 elements which is not the $K$-complex of 
any poset structure. This is because the maximal simplices are $\set{a,b}, \set{a,c}, \set{b,c}$ (and 
all vertices belong to more than one maximal simplex). In general, for the same reasons, 
the boundary of any closed $n$-simplex is not a $K$-complex of any poset structure.

\medskip

In fact, the condition of above is essentially the obstruction to be a $K$-complex:

\begin{thm}\label{associated}
Let $T$ be a finite simplicial complex with vertex set $X$. Then $T$ is the $K$-complex associated to a poset 
structure $\leq$ on $X$  if and only if  for any maximal simplex $s$ of $T$ there exists some
$y\in s$ such that $y\notin s'$ for all maximal simplices $s'\neq s$.
\end{thm}

\begin{proof}
One implication follows immediately from lemma \ref{maximal} and previous remarks. To prove the other implication: Suppose $T$ satisfies
the condition on its maximal simplices. We define a poset structure on $X$ (of length 2) as follows. Let $s_1,\ldots, s_r$ be the maximal simplices of 
$T$. By hypothesis, for each $i=1,\ldots, r$ we can choose some $y_i\in s_i$ such that $y_i\notin s_j$ for all $j\neq i$. These $y_i's$ will be the
maximal elements of the poset. Define the relation $x\leq y$ if $x=y$ or $y=y_i$ for some $i$ and $x\in s_i$.

\medskip

It is not difficult to prove that this is a well defined poset structure on $X$ and that $T$ is the $K$-complex of this stucture. This completes the
proof.
\end{proof}

Of course there is  an analogous result for the associated $L$-complexes.

\begin{remark}
Note that the poset structure constructed in the proof of the Theorem is of length 2. Therefore the $K$-complex of any poset structure on a finite set
$X$ coincides
with the 
$K$-complex of one of length 2. This implies of course that many poset structures on $X$ have the same associated $K$-complex and also that
for some poset structures, the associated $K$-complex does not have the homotopy type of the standard polyhedron $C_X$, since the complex $C_X$ 
 of a poset of length 2  is a graph (=simplicial complex of dimension 1) and any graph has the homotopy type of a bouquet of circles 
 $\bigvee_{\alpha}S^1$.
\end{remark}

\medskip

From Theorem \ref{associated}, one can also deduce:

\smallskip

\begin{remark}
$\k$ is the $K$-complex of $(X,\leq)$ if and only if $(X,\leq)$ has a maximum.
\end{remark}

\bigskip

We investigate now the relationship between the $K$- and $L$-complexes associated to the  relation $\leq$
with the complexes  associated
to the relation $<$.

\smallskip

Let $(X,\leq)$ be a finite poset. As before, we denote by $K$ and $L$ the associated complexes and let  $K'$ and $L'$ be the complexes associated to
the relation 
$<$. It is clear that $K'$ and $L'$ are empty if the poset is discrete, so let us suppose that $X$ is not 
discrete. Moreover, we will suppose for the
next result that no connected component of $X$ is a single point (for example, take $X$ connected with more than one point).

\medskip

By definition, the simplices of $K'$ are the subsets of $X$ consisting of elements $x_0,\ldots,x_n$ such that there is some $y\in X$ with 
$x_i< y$ for all $i$. Clearly $K'< K$. Similarly, $L'<L$.

\bigskip

We will prove that all of them ($K, K', L, L'$) have the same homotopy type. Moreover, we will see that $K$ simplicially 
collapses to $K'$ (which we denote, as 
usual, 
$K\searrow K'$). Similarly one can prove that $L \searrow L'$.

\bigskip

Recall that  
there is an elementary  collapse from a simplicial complex $T$ to a subcomplex $M$ if there is a simplex 
$s$ of $T$ and a vertex $a$ of $T$ not in $s$ such that $T=M\cup as$ and $M\cap as=a\dot{s}$. Equivalently, $T$ elementary collapses to $M$ if
there is a simplex $s$ of $T$ which is contained properly in only one simplex $s'$ and $M=T-\set{s,s'}$.
A finite simplicial complex $K$ collapses to a subcomplex $K'$ if there is a (finite) sequence of elementary collapses from $K$ to $K'$.

\medskip

For a comprehensive exposition on collapses and simple homotopy theory, the classical references are Whitehead's original papers
\cite{Whi, Whi2, Whi4}. The standard references for simple homotopy theory of CW-complexes include Milnor's article \cite{Mil} 
and M.M.Cohen's book \cite{Coh} and for the infinite 
case, Siebenmann's paper \cite{Sie}.

\medskip

In a joint paper with J.A. Barmak \cite{bm2}, we have introduced the notion of collapse for posets  (or equivalently, finite $T_0$-spaces) 
and investigated collapsible posets. An elementary collapse of posets consists of removing just a single point of the poset (which is called 
a \it weak beat point). \rm   Via the associated simplicial complex $C_X$, this notion corresponds to the classical notion of simplicial collapse.

\medskip

To prove our result, we will need the following  basic lemma from \cite{Whi}:

\begin{lemma}[Whitehead \cite {Whi}, Lemma 1]
Let $M, T_0, T$ be finite simplicial complexes such that  $M\cap T_0\subset T$ and $T_0 \searrow T$. Then $M\cup T_0 \searrow M\cup T$.
\end{lemma}

 \begin{thm}
 Let $(X,\leq)$ be a finite poset such that no connected component is a single point. Let $K, L$ be the simplicial complexes associated to $\leq$ and
 $K',L'$ the complexes associated to $<$. Then $K\searrow K'$ and $L\searrow L'$.
 \end{thm}
 
 \begin{proof}
We prove the case $K\searrow K'$, the other case is similar.

\medskip

The simplices of $K-K'$ are exactly the simplices of $K$ containing some maximal element $y$ of $X$. By lemma \ref{maximal}, the maximal simplices of 
$K-K'$ are the simplices of the form $s=\set{y,x_0,\ldots,x_n}$ where $y$ is a maximal element in $X$ and 
the 
set $\set{x_0,\ldots,x_n}$ consists  of all elements of $X$ such that $x_i < y$. Note that the set of elements which are smaller than $y$ is not 
empty by the hypothesis on the connected components of $X$.

\medskip

Since all the faces of the maximal simplex $s=\set{y,x_0,\ldots,x_n}$ which contain the vertex $y$ are not faces of any other maximal simplex $s'$, 
we can suppose without loss of generality, that $K-K'$ contains only one maximal simplex $s$.

\medskip

Let $T_0$ be the closed $(n+1)$-simplex spanned by $s=\set{y,x_0,\ldots,x_n}$ and let $T$ be the closed $n$-simplex spanned by $\set{x_0,\ldots,x_n}$.
Since $T_0\searrow T$ (\cite{Whi}, Lemma 2), by the lemma of above we have
$$K=K'\cup T_0 \searrow K'\cup T=K'$$
\end{proof}

We finish the paper with a result that relates the $K$-complexes of two posets $X, Y$ with  \it closed \rm relations
$R\subset X\times Y$. Given two posets $X, Y,$ a relation $R\subset X\times Y$ is called closed if it satisfies the following property: For 
any $(x,y)\leq (x',y')\in X\times Y$, if $(x,y)\in R$, then $(x',y')\in R$ (in the finite case, this is equivalent to $R$ being
 a closed subset of the finite space $X\times Y$ with the
product topology).

There is a well-known result by Quillen \cite{Qui} (see also \cite{Bjo0}, Thm. 10.10) which relates the closed relations with the simplicial complexes
of finite chains $C_X$ and $C_Y$. Explicitly,

\begin{thm}[Quillen]
Let $X$ and $Y$ be posets and let $R\subset X\times Y$ be a closed relation. For any  $x\in X$ and $y\in Y$, let $C_x$ be the simplicial complex
of nonempty finite chains of the subposet
 $S_x=\set{y\in Y,\ xRy}\subset Y$ and $C_y$ the corresponding simplicial complex of the subposet $S_y=\set{x\in X,\ xRy}\subset X$. 
 If $C_x$ and $C_y$ are contractible for all $x$ and $y$, then $C_X$ and $C_Y$ are homotopy equivalent.
 \end{thm}
 
We want a similar result for the $K$-complexes. Unfortunately the analogous result (replacing $C_x, C_y, C_X$ and $C_Y$ by the corresponding
$K$-complexes) is not valid, as we show in the following example.

\begin{example}
Consider the following posets $X$ and $Y$
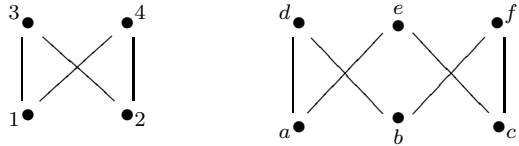
\begin{displaymath}
 \xymatrix@C=10pt{ ^3 \bullet \ar@{-}[d] \ar@{-}[drr] & & \bullet ^4 \ar@{-}[lld]  \ar@{-}[d]  \\
		_1 \bullet & & \bullet _2  } 
\qquad \qquad
 \xymatrix@C=10pt{ ^d \bullet \ar@{-}[d] \ar@{-}[drr] &  & {\overset{e}{\bullet}} \ar@{-}[dll]\ar@{-}[drr] & & \bullet ^f \ar@{-}[lld]  \ar@{-}[d]  \\
		_a \bullet & & {\underset{b}{\bullet}} & & \bullet _c  } 
\end{displaymath}

and the closed relation 
$$R=\set{(1,d),(2,e),(3,b),(3,c),(3,d),(3,e),(3,f),(4,a),(4,d),(4,e)}$$
It is easy to see that the $K$-complex of each $S_x$ and $S_y$ is contractible but the $K$-complexes of $X$ and $Y$ are not homotopy equivalent:
The first one is contractible and the second one is homotopy equivalent to $S^1$. This is also an example of two posets with $C_X \simeq C_Y$ but
with $K$-complexes of different homotopy types.
 \end{example}

However we obtain the following weak version of the Theorem for the $K$-complexes (for the $L$-complexes one has of course a similar result).

\begin{thm}
Let $X$ and $Y$ be posets and let $R\subset X\times Y$ be a closed relation. If the subposets $S_x$ and $S_y$ have maximum for all $x\in X$ and 
$y\in Y$, then the $K$-complexes $K_X$ and $K_Y$ are homotopy equivalent. 
\end{thm}

Before we proceed with the proof we make a couple of remarks. Note first that the simplicial complex $K_Y$ is the $K$-complex of the poset $(Y,\leq)$ (it
is not the $K$-complex of the relation $(Y,R)$ as in the first part of the paper).

\medskip

Note also that the hypothesis required for this Theorem is strictly stronger than the hypothesis of Quillen's result: If $S_x$ has a maximum, then 
in particular
$C_x$ and $K_x$ are contractible.

\medskip

\begin{proof}
Consider $R$ as a poset and denote by $K_R$ the $K$-complex of $R$. We will prove that the proyection $R\to X$ (resp. $R\to Y$) induces a
 homotopy
equivalence $K_R\to K_X$ (resp. $K_R\to K_Y$). 

\bigskip

By Quillen's Theorem A \cite{Qui0}, it suffices to prove that $p^{-1}(s)\subset K_R$ is contractible for every closed simplex $s$ of $K_X$.

\medskip

Let $s=\set{x_0,\ldots,x_n}$ be a simplex of $K_X$, then there exists $x'\in X$ such that $x_i\leq x'$ for all $i$. Let $(x,y)$ be a vertex of
$p^{-1}(s)\subset K_R$. Therefore $x=x_i$ for some $i=0,\ldots,n$ and $xRy$. Since the relation is closed, we have that $x'Ry$. 
Denote by $y'$ the maximum
element of $S_{x'}$. Then $(x,y)\leq (x',y')$.

\medskip

We have proved that all the vertices of $p^{-1}(s)\subset K_R$ are smaller than $(x',y')$. By definition of the $K$-complex, this implies that
$p^{-1}(s)\subset K_R$ is a closed simplex (or a generalized simplex in the infinite case) and therefore contractible.
\end{proof}

\email{ gminian$@$dm.uba.ar}
\end{document}